%&bigtex

%%   Amstex document
%%   Updated 1/26/04 by SF:  Corrected typos, included lemma
%%   about the determinant for the chain case.
%%
%%   Updated 1/28/04 by JL:  Fixed typographical issues,
%%   numbering of lemmas and equations, rewrote proof of
%%   Theorem 1 using results from new reference [M]
%%
%%   Update 2/23/04 by CB: included example with diagrames
%%   for the boolean algebra of rank 2
%%
%%   Update 3/12/04 by CB: corrections following discussion on 3/11/04
%%
%%   Update 3/15/04 by CB: intro to section 2, more corrections

\input amstex
\documentstyle{amsppt}
\input xy
\xyoption{all}
\xyoption{graph}
\xyoption{tips}
\magnification=1200
\hcorrection{.25in} \advance\vsize-.75in
\def\proof{{\smc Proof. }}
\def\endproof{$\quad \square$}

\overfullrule=0pt
\topmatter
\title{Determinants Associated to Zeta Matrices of Posets}\endtitle
\leftheadtext{Determinants of Zeta Matrices} \rightheadtext{Determinants of Zeta
Matrices} \keywords poset, zeta function, M\"obius function\endkeywords \subjclass
15A15, 05C20, 05C50,  06A11  \endsubjclass
\author
Cristina M. Ballantine\\
Sharon M.  Frechette\\
John B. Little\\
\\
College of the Holy Cross\\
\endauthor
\address Department of Mathematics and Computer Science,
College of the Holy Cross,
Worcester, MA 01610\endaddress
\email (cballant,sfrechet,little)\@mathcs.holycross.edu \endemail
\date  March 15, 2004 \enddate

\abstract We consider the matrix ${\frak Z}_P=Z_P+Z_P^t$, where the
entries of
$Z_P$ are the values of the zeta function of the
finite poset $P$.
We give a combinatorial interpretation of the
determinant of ${\frak Z}_P$
and establish a recursive formula for this
determinant in the case
in which $P$ is a boolean algebra.
\endabstract
\endtopmatter

\document
\SelectTips{cm}{} \objectmargin={1pt}

 \heading \S 1. Introduction \endheading
\bigskip

The theory of partially ordered sets (posets) plays an important
role
in enumerative combinatorics and the M\"obius inversion
formula for posets
generalizes several fundamental theorems
including the number-theoretic
M\"obius inversion theorem. For a
detailed review of posets and M\"obius
inversion we refer the
reader to \cite{S1}, chapter $3$, and \cite{Sa}.
Below we provide
a short exposition of the basic facts on the subject
following
\cite{S1}.
\bigskip
A partially ordered set ({\it poset}) $P$
is a set which, by abuse
of notation, we also call $P$ together with a
binary relation,
called a {\it partial order} and denoted $\leq$, satisfying:
\roster
\item $\ x \leq x$ for all $ x \in P$ (reflexivity).

\item If $x\leq y$ and $y \leq x$, then $x=y$ (antisymmetry).

\item If $x\leq y$ and $y \leq z$, then $x \leq z$ (transitivity).
\endroster
Two elements $x$ and $y$ are {\it comparable} if $x \leq y$ or $y
\leq x$.
Otherwise they are {\it incomparable}. We write $x<y$ to
mean $x \leq y$ and
$x \neq y$.

\bigskip

\noindent
{\bf Examples.}
\roster
\item Let $n \in \Bbb{N}$.
The set
$[n]=\{1,2,\ldots, n\}$ with the usual order forms a poset in
which
any two elements are comparable. Such a poset is called a
{\it chain} or a
totally ordered set.

\item Let $n \in \Bbb{N}$. Consider the poset
$P_n$ of subsets of
$[n]$ under the inclusion relation. This poset is called a
{\it
boolean algebra} of rank $n$. In \cite{S1} it is denoted by
$2^{[n]}$.

\item Let $n \in \Bbb{N}$. The set $D_n$ of all positive divisors of
$n$ forms
a poset under the order defined by $i \leq j$ in $D_n$
if $j$ is divisible by
$i$. This poset is called the {\it divisor
poset}.

\endroster

A closed {\it interval} $[x,y]$ is defined whenever $x \leq y$ by
$[x,y]=\{z \in P: x \leq z \leq y\}$. The empty set is not regarded
as an
interval. Denote by $Int(P)$ the set of intervals of $P$. If
$f: Int(P) \rightarrow \Bbb{C}$, we write $f(x,y)$ for $f([x,y])$.
\bigskip
The {\it incidence algebra} $I(P)$ of $P$ is the $\Bbb{C}$-algebra
of all
functions $f: Int(P) \rightarrow \Bbb{C}$ under the
convolution given by
$$fg(x,y)= \sum_{x \leq z \leq
y}f(x,z)g(z,y).$$

It is known that $I(P)$ is an associative $\Bbb{C}$-algebra with identity
$$\delta(x,y)= \cases 1 & \text{if }  x=y,\cr
                      0 & \text{if } x \neq y.\cr \endcases$$

The {\it zeta function} $\zeta \in I(P) $ of a poset $P$ is defined by
$\zeta(x,y)=1$ for all  $x \leq y$ in $P$. If $P$ is a locally finite poset (i.e.
every interval in $P$ is finite), the zeta function $\zeta$ is invertible in the
algebra $I(P)$. Its inverse is called the {\it M\"obius function} of $P$ and is
denoted by $\mu$. Note that one can define $\mu$ inductively by
\bigskip

$\mu(x,x)=1$, for all $x \in P,$
\bigskip

${\displaystyle \mu(x,y)=-\sum_{x \leq z < y}\mu(x,z)}$, for all $x<y$ in $P$.
\bigskip
\proclaim{Proposition 1}(M\"obius Inversion Formula) Let $P$ be a
finite
poset and
\newline
$f,g: P \rightarrow \Bbb{C}$. Then
$$g(x)=\sum_{y \leq x}f(y) \text{ for all } x \in P
\Leftrightarrow
f(x)=\sum_{y \leq x} g(y) \mu(y,x) \text{ for all } x \in P.$$
\endproclaim
\bigskip

\noindent
\proof See \cite{S}, pg. 116.
\endproof

\bigskip

There is also a dual form of the M\"obius Inversion Formula.
Moreover, the formula works for more general posets. All that is needed
is
that every principal order ideal be finite.

\bigskip

\noindent {\bf Remark.} If $P$ is the divisor poset of Example (3), the
proposition above becomes the number-theoretic M\"obius Inversion
Theorem.
If $P$ is the boolean algebra of Example (2), the
proposition gives the
Principle of Inclusion-Exclusion. Finally,
if $P$ is the chain of Example
(1), the proposition becomes the
Fundamental Theorem of the Difference
Calculus. For the
aforementioned theorems see \cite{Sa}.

\bigskip

For the remainder of the  article, $P$ will be a poset with $n$
elements
and the partial order denoted by $\leq$.  We choose a
labelling
$x_1, x_2, \ldots, x_n$ of the elements of $P$ such that
$x_i < x_j \Longrightarrow i < j$.

\bigskip
\noindent {\bf Definition.} The {\it zeta matrix} $Z_P$ of a poset
$P$
is defined as the $n \times n$ matrix with entries
$$(Z_P)_{ij}= \cases 1 & \text{if }  x_i \leq x_j\cr
                     0 & \text{otherwise} \cr\endcases$$
\bigskip
Observe that, with the chosen labelling, the zeta matrix is unipotent upper
triangular. Its non-zero entries are the values of the zeta function.

\bigskip
We define the matrix $\frak{Z}_P$ by $\frak{Z}_P=Z_P + Z_P^t$.
In
Section \S 2 we give a combinatorial interpretation of the determinant
of $\frak{Z}_P$. Section \S 3 refines the interpretation in the case
in which
$P$ is a chain. The value of the determinant in this case is
$n+1$. The main
theorem of the paper evaluates the determinant of
$\frak{Z}_n := \frak{Z}_{P_n}$ when $P_n$ is the boolean algebra of
rank $n$.
More specifically, in Section \S 4, we prove the following
recursive formula
on $n$.

\bigskip
\proclaim{Main Theorem}  If $n \ge 3$ is odd, then $\det(\frak{Z}_n) = 0$. If $n$ is
even, then $\det(\frak{Z}_n) = 2^{\alpha_n}$ where $\alpha_2 = 2$, and $\alpha_n = 4
\alpha_{n-2} - 2$ for $n \ge 4$ and even.
\endproclaim
\bigskip

 Consider also the matrix $\frak{M}_P$ defined by
$\frak{M}_P=M_P + M_P^t$, where $M_P=Z_P^{-1}$.  The non-zero entries of $M_P$ are
the values of the M\"obius function.  We refer to $M_P$ as the {\it M\"obius matrix}
of the poset $P$.  We have the following theorem.
\bigskip

\proclaim{Theorem 1} $\det(\frak{M}_P)=\det(\frak{Z}_P)$.
\endproclaim
\bigskip

\noindent \proof  The theorem is a direct consequence of the following lemma.
\endproof

\proclaim{Lemma 1} Let $U$ be an $n \times n$ matrix such that $\det(U)=1$ and let
$V=U^{-1}$. Then, $\det(U+U^t) =\det(V+V^t)$.
\endproclaim
\bigskip

\noindent
\proof We have
$$V^t+V=(U^{-1})^t+U^{-1}=(U^{-1})^t U
U^{-1}+(U^t)^{-1} U^t U^{-1} =(U^{-1})^t[U+U^t]U^{-1}.$$ Thus
$\det(V+V^t)=\det((U^{-1})^t)\det(U+U^t)\det(U^{-1})$ and since $\det(U)=1$, we have
$\det(V+V^t)=\det(U+U^t)$.
\endproof

\bigskip

\bigskip
\heading \S 2. Combinatorial interpretations of $\det(\frak{Z}_P)$ \endheading

\bigskip

In this section, we will give two related combinatorial interpretations of
$\det(\frak{Z}_P)$. The first involves collections of disjoint cycles in a directed
graph $D_P$ associated to $P$. The second involves cycle decompositions of elements
of a certain subgroup $S_n^P$ of the symmetric group $S_n$ associated to $P$.
\bigskip

Consider a poset $P$ as in the previous section with $|P| = n$. The matrix $Y_P =
Z_P - I_n$, in which the diagonal entries of $Z_P$ are replaced by $0$, can be
interpreted as the adjacency matrix of a {\it directed graph} (digraph) $G_P$
associated to the strict order relation $x < y$ in $P$.    The vertices of $G_P$ are
the elements of $P$, and there is a directed edge from $x$ to $y$ if and only if $x
< y$. Similarly, the matrix ${\frak Y}_P = Y_P + Y_P^t$ is the adjacency matrix of
the directed graph $D_P$ in which there are edges in {\it both} directions between
each pair of distinct comparable elements $x,y \in P$. Then we have
$$\frak{Z}_P = Z_P + Z_P^t = Y_P + Y_P^t + 2I_n = {\frak Y}_P + 2I_n,$$
and $\det(\frak{Z}_P) = \det({\frak Y}_P + 2I_n) = \chi(-2)$ where $\chi(t)$ is the
characteristic polynomial $\chi(t) = \det({\frak Y}_P - tI_n)$ of the matrix ${\frak
Y}_P$.  As in \cite{M}, we will also call this the characteristic polynomial of the
graph $D_P$. \bigskip

If $P=\{x_1, \ldots, x_n\}$ and $\frak{Y}_P-tI_n=(m_{ij})$, then by the definition
of the determinant

$$\chi(t)=\sum_{\sigma \in S_n} (-1)^{sign(\sigma)} \prod_{j=1}^n m_{j \sigma(j)}.$$

\bigskip

Then we have
$$\prod_{j=1}^n m_{j \sigma(j)} =  \cases 0 & \text{if }  x_j \text{ not comparable to }
x_{\sigma(j)} \text{ for some } j,\cr
                      (-t)^{k_{\sigma}} & \text{ otherwise},
                      \cr \endcases$$ where $k_{\sigma}$ is the number of fixed points
of the permutation $\sigma$.

\bigskip

Let $S_n^P=\{\sigma \in S_n \text{ : } x_j \text{ comparable in } P \text{ to }
x_{\sigma(j)},\, \text{ for all } j=1, \ldots , n\}$. Then combining the last two
displayed formulas $$\chi(t)= \sum_{\sigma \in S_n^P}(-1)^{sign(\sigma)}
(-t)^{k_\sigma}.\leqno(1)$$

As shown in the work of Harary and Mowshowitz (see \cite{H} and \cite{M}),
there is an interesting relation between the coefficients of $\chi(t)$
and cycles in the graph $D_P$.  Namely, if we write
$$\chi(t) = \sum_{i=0}^n  a_i (-t)^{n-i},$$
then $a_0 = 1$ (the coefficient of $(-t)^n$ equals $1$ since the identity
permutation in $S_n$ is the only permutation with exactly $n$ fixed points) and
Theorem 1 of \cite{M} demonstrates that for $i \ge 1$,
$$a_i = \sum \left(\prod_{j=1}^r (-1)^{i_j + 1}\right) f_{D_P}(i_1,\ldots,i_r),\leqno(2)$$
where the summation is taken over all partitions $i = i_1 + \cdots + i_r$ with $i_j
> 0$ for all $j$, and $f_{D_P}(i_1,\ldots,i_r)$ is the number of collections of disjoint
directed cycles in $D_P$ of lengths $i_1,\ldots,i_r$.  Thus

$$a_i = \sum  (-1)^{i+r} f_{D_P}(i_1,\ldots,i_r).$$

Note that in fact $i_j \ge 2$ for all $j$ since $D_P$ contains no loops (no edges
from a vertex to itself). Therefore $a_1=0$, which corresponds to the fact that
there are no permutations in $S_n$ with exactly $n-1$ fixed points.

\bigskip
\noindent {\bf Example.} To illustrate the relation between $(1)$ and $(2)$, let
$P=P_2$ be the boolean algebra of rank $2$. The ordering of the elements of $P$ is
given by the correspondence $$ \emptyset \longleftrightarrow 00, \ \ \{1\}
\longleftrightarrow 01, \ \ \{2\} \longleftrightarrow 10, \ \ \{1,2\}
\longleftrightarrow 11.$$ Thus we have the labelling
$$P=\{x_1=00, x_2=01, x_3=10, x_4=11\}.$$  Then $$\frak{Z}_P= \pmatrix 2 & 1 & 1 & 1
\cr 1 & 2 & 0 & 1 \cr 1 &  0 & 2 & 1 \cr 1 & 1 & 1 &  2 \cr \endpmatrix$$ and
$$\chi(t)=\det(\frak{Y}_P-tI_4)=\det\pmatrix -t & 1 & 1 & 1
\cr  1 & -t & 0 & 1 \cr 1 &  0 & -t & 1 \cr 1 & 1 & 1 &  -t \cr \endpmatrix$$ $$=
\sum_{i=0}^4 a_i (-t)^{4-i}= (-t)^4-5(-t)^2+4(-t),$$ where the $a_i$'s are defined
as above, so $a_0=1, \ a_1=0, \ a_2=-5, \ a_3=4, \ a_4=0$.
\bigskip

The digraphs $G_P$ and $D_P$ are given below.
\bigskip

$G_P$: \ \ \ \ \xygraph{ !{0;(.77,-.77):0}*+++!D{11}
*{\bullet}"11"(:@{-}|@{<}[r(1.3)]*++!L{10}*{\bullet}"10":@{-}|@{<}[d(1.3)]
*++!U{00}*{\bullet}"00"
:@{-}|@{>}[l(1.3)]*++!R{01}*{\bullet}"01":@{-}|@{>}[u(1.3)]"11":@{-}|@{<}[dr(1.3)]"00")}\
\ \ \ \ \ \ \ \ $D_P$: \ \ \ \ \xygraph{ !{0;(.77,-.77):0}!~:{@{-}|@{>}}*+++!D{11}
*{\bullet}"11"(:@/^/[r(1.3)]*++!L{10}*{\bullet}"10":@/^/[d(1.3)]*++!U{00}*{\bullet}"00"
:@/^/[l(1.3)]*++!R{01}*{\bullet}"01":@/^/[u(1.3)]"11") :@/^/[dr(1.3)]"00"
:@/^/[u(1.3)]"10":@/^/[l(1.3)]"11":@/^/[d(1.3)]"01":@/^/[r(1.3)]"00":@/^/[ul(1.3)]"11"}
\bigskip

We  obtain one term of $(-t)^2$ for each permutation in $S_4^P$ with exactly $2$
fixed points. Since each part of the partition must be at least $2$, the only
allowable partition of $2$ is $2=2$ and ${\displaystyle f_{D_P}(2)}$ counts the
following collections of cycles in $D_P$ which correspond to the transpositions in
$S_4^P$:
\bigskip

$\ $ \xygraph{ !{0;(.77,-.77):0}!~:{@{-}|@{>}}*+++!D{11}
 *{\bullet}"11"{}[r]*++!L{10}*{\bullet}"10"
{}[d]*++!U{00}*{\bullet}"00" :@/^/[l]*++!R{01}*{\bullet}"01":@/^/[r]"00"} \ \ $\ \ \
\ \ \ $ \ \ \ \xygraph{ !{0;(.77,-.77):0}!~:{@{-}|@{>}}*+++!D{11}
*{\bullet}"11"{}[r]*++!L{10}*{\bullet}"10" {}[d]*++!U{00}*{\bullet}"00"
{}[l]*++!R{01}*{\bullet}"01":@/^/[u]"11":@/^/[d]"01"} \ \ \ \ \ \ $\ \ \ \ $
 \xygraph{
!{0;(.77,-.77):0}!~:{@{-}|@{>}}*+++!D{11}
*{\bullet}"11"{}[dr]*++!U{00}*{\bullet}"00"{}[l]*++!R{01}*{\bullet}"01" {}[u]
"11"{}[r]*++!L{10}*{\bullet}"10" :@/^/[d]"00":@/^/[u]"10" }

\bigskip

\bigskip

$\ $  \xygraph{ !{0;(.77,-.77):0}!~:{@{-}|@{>}}*+++!D{11}
*{\bullet}"11":@/^/[r]*++!L{10}*{\bullet}"10":@/^/[l]"11"{}[d]
*++!R{01}*{\bullet}"01" {}[r]*++!U{00}*{\bullet}"00" }\ \ \ \ $\ \ \ \ \ \ \ \ \ \ \
\ \ $\ \ \xygraph{ !{0;(.77,-.77):0}!~:{@{-}|@{>}}*+++!D{11}
*{\bullet}"11":@/^/[dr]*++!U{00}*{\bullet}"00":@/^/[ul]"11"
{}[r]*++!L{10}*{\bullet}"10"{}[dl] *++!R{01}*{\bullet}"01" }
\bigskip

This corresponds to the coefficient $a_2=-5$ .
\bigskip

We obtain one term of $(-t)^1$ for each permutation in $S_4^P$ with exactly $1$
fixed point. The only allowable partition of $3$ is $3=3$ and ${\displaystyle
f_{D_P}(3)}$ counts the following collections of cycles in $D_P$ which correspond to
the $3$-cycles in $S_4^P$:
\bigskip

$\ $ \xygraph{ !{0;(.77,-.77):0}*+++!D{11}
*{\bullet}"11"(:@{-}|@{>}[d]*++!R{01}*{\bullet}"01":@{-}|@{>}[r]*++!U{00}
*{\bullet}"00" :@{-}|@{>}[ul]"11":@{}[r]*++!L{10}*{\bullet}"10")} \ \ \ \ $\ \ \ \
$\ \ \  \xygraph{ !{0;(.77,-.77):0}*+++!D{11}
*{\bullet}"11"(:@{-}|@{>}[dr]*++!U{00}*{\bullet}"00":@{-}|@{>}[l]
*++!R{01}*{\bullet}"01" :@{-}|@{>}[u]"11":@{}[r]*++!L{10}*{\bullet}"10")}
\bigskip

\bigskip

$ \ $ \xygraph{ !{0;(.77,-.77):0}*+++!D{11}
*{\bullet}"11"(:@{-}|@{>}[r]*++!L{10}*{\bullet}"10":@{-}|@{>}[d]
*++!U{00}*{\bullet}"00" :@{-}|@{>}[ul]"11":@{}[d]*++!R{01}*{\bullet}"01")} \ \ \ \ $
\ \ \  $ \ \ \ \xygraph{ !{0;(.77,-.77):0}*+++!D{11}
*{\bullet}"11"(:@{-}|@{>}[dr]*++!U{00}*{\bullet}"00":@{-}|@{>}[u]
*++!L{10}*{\bullet}"10" :@{-}|@{>}[l]"11":@{}[d]*++!R{01}*{\bullet}"10")}

\bigskip

This corresponds to the coefficient $a_3=4$.\bigskip

Finally, we  obtain one  term of $(-t)^0=1$ for each permutation in $S_4^P$ with no
fixed points. Such permutations are either $4$-cycles or products of $2$
transpositions.  The only allowable partitions are $4=4$ and $4=2+2$. Then
${\displaystyle f_{D_P}(4)}$ counts the following collections of cycles in $D_P$
(with a $-$ sign) which correspond to the $4$-cycles in $S_4^P$:
\bigskip

$\ $ \xygraph{ !{0;(.77,-.77):0}!~:{@{-}|@{>}}*+++!D{11}
*{\bullet}"11"(:@{-}[d]*++!R{01}*{\bullet}"01":@{-}[r]*++!U{00}*{\bullet}"00"
:@{-}[u]*++!L{10}*{\bullet}"10":@{-}[l]"11")}\ \ \ \ $\ \ \ \ $\ \ \
 \xygraph{
!{0;(.77,-.77):0}!~:{@{-}|@{>}}*+++!D{11}
*{\bullet}"11"(:@{-}[r]*++!L{10}*{\bullet}"10":@{-}[d]*++!U{00}*{\bullet}"00"
:@{-}[l]*++!R{01}*{\bullet}"01":@{-}[u]"11")}
\bigskip

On the other hand,  ${\displaystyle f_{D_P}(2,2)}$ counts the following collections
of cycles in $D_P$ (with a $+$ sign) which correspond to the permutations in $S_4^P$
which are products of $2$ transpositions:
\bigskip

$\ $ \xygraph{ !{0;(.77,-.77):0}!~:{@{-}|@{>}}*+++!D{11}
*{\bullet}"11"(:@/^/[r]*++!L{10}*{\bullet}"10":@/^/[l]"11"
{}[dr]*++!U{00}*{\bullet}"00" :@/^/[l]*++!R{01}*{\bullet}"01":@/^/[r]"00")}\ \ \ \ \
$ \ \ \ $\ \ \  \xygraph{ !{0;(.77,-.77):0}!~:{@{-}|@{>}}*+++!D{11}
*{\bullet}"11"(:@/^/[d]*++!R{01}*{\bullet}"01":@/^/[u]"11"
{}[dr]*++!U{00}*{\bullet}"00" :@/^/[u]*++!L{10}*{\bullet}"10":@/^/[d]"00")}\bigskip

Thus $a_4=-2+2=0$. \bigskip

Returning to our general situation, substituting from (1) we have
$$\eqalign{
\chi(-2) &= \sum_{i=0}^n a_i 2^{n-i}\cr
         &= 2^n + \sum_{i=1}^n \left(\sum_{\sum i_j = i} (-1)^{i+r} f_{D_P}(i_1,\ldots,i_r)
         \right) 2^{n-i}.\cr}$$
For each collection of $r$ disjoint directed cycles there are $n-i$ vertices not
contained in the cycles. Let $s = n - i + r$ be the total number of connected
components of the union of the cycles and the disconnected vertices.  We  reindex
the last sum as a sum over $s$:
$$\eqalign{
\chi(-2) &= \sum_{s=1}^n (-1)^{n-s} \left(\sum_{r,(i_j) \atop s = r + n - \sum i_j}
f_{D_P}(i_1,\ldots,i_r)2^{s-r}\right)\cr &=\sum_{s=1}^n (-1)^{n-s} c_s.\cr}$$ (The
sign comes from the fact that $i + r = s - n + 2i$, so $(-1)^{i+r} = (-1)^{n-s}$.)
The coefficient
$$c_s = \sum_{r,(i_j)\atop s = r + n - \sum i_j} f_{D_P}(i_1,\ldots,i_r)2^{s-r}$$
can also be interpreted as the integer
$$\sum_{r,(i_j)\atop s = r + n - \sum i_j} f_{\overline{D_P}}(i_1,\ldots,i_r,1,\ldots,1),$$
(where there are $s-r$ ones) counting the number of
collections of disjoint directed cycles
of lengths $i_1,\ldots,i_r,1,\ldots,1$ containing {\it all}
the vertices in a directed graph
$\overline{D_P}$ obtained from $D_P$ by adding two loops at
each vertex (corresponding to the diagonal entries $2$ in the
matrix ${\frak Z}_P = Z_P + Z_P^t$).  Hence $c_s$ is
the {\it number of collections of disjoint directed cycles in the graph
$\overline{D_P}$ having precisely $s$ connected components and containing
all the vertices.}  Note that
$s - r$ is the number of loops contained in each such collection.

To summarize, we have proven the following theorem.

\bigskip
\proclaim{Theorem 2}  With the above notation, for all finite partially ordered sets
$P$,
$$\det({\frak Z}_P) = \sum_{s=1}^n (-1)^{n-s} c_s.$$
\endproclaim

In terms of permutations in $S_n^P$, if  $1 \leq s \leq n-1$, we have

$$ c_s = \sum_{\ell=0}^{s-1} \left|\left\{\sigma \in S_n^P :  \matrix \sigma \text{ is a
product of $s-\ell$ disjoint }\cr \text{cycles of length } > 1
\endmatrix \right\} \right|2^\ell \leqno(3)$$
and $c_n=2^n$.

\bigskip
\heading \S 3. The chain case \endheading

\bigskip
Consider the special case in which  $P$ is the chain of length
$n$.
Then $P=[n]=\{1,2, \ldots, n\}$ with the usual ordering. All
elements
of $P$ are comparable, and we have
$$Z_{[n]} = \pmatrix 1 & 1 & \cdots & 1\cr
                       & 1 & \cdots & 1\cr
                       &   & \ddots & \vdots\cr
                       &   &        & 1\cr\endpmatrix
\qquad
\text{and}
\qquad
M_{[n]} = Z_{[n]}^{-1} = \pmatrix 1 & -1 &  0 & \cdots & 0\cr
                                           &  1 & -1 & \cdots & 0\cr
                                           &  & \ddots & \ddots & \vdots\cr
                                           &  &  &  1 & -1\cr
                                           &  &  &  & 1\cr\endpmatrix,$$
where all lower-triangular entries in each matrix are zero.
\bigskip

\proclaim{Lemma 2}  With the above notation
$$\det(\frak{Z}_{[n]}) = \det(\frak{M}_{[n]}) = n+1.$$
\endproclaim

\bigskip
\noindent \proof The first equality holds by Theorem 1. We prove the second equality
by induction on $n$. We have $\frak{M}_{[1]}=(2)$, hence the lemma holds when $n=1$.
Suppose the lemma holds for all $k < n$, and consider
$$\frak{M}_{[n]} = \pmatrix
 2 & -1 &  & & \cr
                            -1&  2 & -1 & & \cr
                            &  -1 & \ddots & \ddots & & \cr
                            &  & \ddots  & 2 & -1\cr
                            &  &  & -1  & 2\cr\endpmatrix,$$
where all remaining entries are zero.  Expanding along the first row, we have
$\det(\frak{M}_{[n]}) = 2 \det(\frak{M}_{[n-1]})-(-1)\det(B_n)$ where
$$B_n = \pmatrix -1 & -1 & 0 &\cdots & 0\cr
                  0 &  &  & & \cr
                  \vdots&  & \frak{M}_{[n-2]}  & & \cr
                  0 &  &  &   & \cr\endpmatrix.$$
Now expanding along the first column in $B_n$, we have $\det(B_n) = -
\det(\frak{M}_{[n-2]})$.  Thus we have
$$\det(\frak{M}_{[n]}) = 2
\det(\frak{M}_{[n-1]}) - \det(\frak{M}_{[n-2]}) = 2n - (n-1) = n+1$$
by
applying the inductive hypothesis.
\endproof
\bigskip

All elements of $P$ are comparable in the chain case, and so the digraph
$\overline{D_P}$ is complete with double edges between all vertices and double loops
at all vertices. Considering permutations written as a product of disjoint cycles,
denote by $\Gamma(n; \gamma_1, \gamma_2, \ldots, \gamma_n)$ the number of
permutations of $[n]$ with exactly $\gamma_i$ cycles of length $i$ for each $1 \leq
i \leq n$. By \cite{S1}, Proposition 1.3.2 we have
$$  \Gamma(n; \gamma_1, \gamma_2, \ldots, \gamma_n)=
\frac{n!}{1^{^{\gamma_1}} \gamma_1! \ 2^{^{\gamma_2}} \gamma_2!\ \ldots \
n^{^{\gamma_n}} \gamma_n!}.$$
\bigskip
Then, in the chain case, the formula for $c_i$ in (3) becomes
$$c_i = \sum_{{(\gamma_1, \ldots,
\gamma_n)}\atop {{\sum j \gamma_j=n}\atop {\sum
\gamma_j=i}}}\Gamma(n; \gamma_1, \gamma_2, \ldots,
\gamma_n)2^{\gamma_1}.$$

\bigskip
The $c_i$'s may be interpreted as {\it modified} Stirling numbers of the second
kind. Recall that the Stirling numbers of the second kind $\left[\matrix n\cr m\cr
\endmatrix\right]$ count the number of permutations of $n$ elements with $m$
disjoint cycles and are given by
$$\left[\matrix n \cr m\cr \endmatrix\right]=
\sum_{{\sum m_i=m}\atop {\sum i m_i = n}}\frac {n!}{1^{m_1}
m_1! \ 2^{m_2} m_2! \ \ldots \ k^{m_k} m_k!} =
\sum_{{\sum m_i=m}\atop{\sum i m_i = n}}\Gamma(n; m_1, \ldots,
m_k). $$
\bigskip

\proclaim{Corollary 1} The alternating sum of the modified
Stirling numbers
satisfies the following.

$$\sum_{i=1}^n (-1)^{n-i}
\sum_{{\sum j \gamma_j=n}\atop{\sum \gamma_j=i}}\Gamma(n;
\gamma_1, \ldots , \gamma_n)2^{\gamma_1}=n+1$$
\endproclaim

\heading
\S 4.  The Boolean algebra case \endheading

\bigskip
Let $[n] = \{1,2, \ldots, n\}$ and consider the poset $P_n=2^{[n]}$ of subsets of
$[n]$ under the inclusion relation.   Let $\frak{Z}_n $ be the matrix defined in \S
1. The main result of this section if the following theorem.
\bigskip
\proclaim{Theorem 3}  If $n \ge 3$ is odd, then $\det(\frak{Z}_n) = 0$. If $n$ is
even, then $\det(\frak{Z}_n) = 2^{\alpha_n}$ where $\alpha_2 = 2$, and $\alpha_n = 4
\alpha_{n-2} - 2$ for $n \ge 4$ and even.
\endproclaim

\bigskip
\noindent \proof The proof will rely on  several lemmas.  First, we identify a
particularly useful labelling of $P_n = 2^{[n]}$ for our purposes, and we consider
only this labelling in the following. Each subset $A \in P_n$ will be encoded as a
binary vector $v(A)$ of length $n$:
$$v(A) =  (v_n,v_{n-1},\ldots,v_1)$$
where
$$v_i = \cases 1 & \text{if } i \in A,\cr
               0 & \text{if not.}\cr\endcases$$
Our labelling of $P_n$ induces the usual numerical ordering when we interpret $v(A)$
as the binary expansion of an integer $m$, with $0 \le m \le 2^n - 1$.

Using this labelling yields an interesting recursive structure in the matrices
$\frak{Z}_n$.

\bigskip
\proclaim{Lemma 3}  The matrices $Z_n$ and $\frak{Z}_n$ have the following
properties. \roster \item The  entries of $Z_n$ above the diagonal are the first
$2^n$ rows in the Pascal triangle modulo 2. \item  For $n \ge 2$, $Z_n$ and
$\frak{Z}_n$ have  block decompositions:
$$Z_n = \pmatrix Z_{n-1} & Z_{n-1}\cr
                          0& Z_{n-1}\cr \endpmatrix \ \text{   and   }\
                          \frak{Z}_n = \pmatrix \frak{Z}_{n-1} & Z_{n-1}\cr
                          Z_{n-1}^t & \frak{Z}_{n-1}\cr \endpmatrix.$$(This statement also
holds with $n = 1$ if we take $\frak{Z}_0 = 2$, $Z_0 = 1$.) \item  The $Z_n$ matrix
sequence can be generated by a recursive procedure as follows.  Given $Z_{n-1}$, to
form $Z_n$ we replace each entry $1$ by a $2\times 2$ block $\pmatrix 1&1\cr
0&1\cr\endpmatrix$ and each entry $0$ by a $2\times 2$ zero matrix. From part $(1)$,
we get a similar recursive procedure for the $\frak{Z}_n$ sequence.
\endroster
\endproclaim

\bigskip
\noindent \proof These properties  follow directly from the definition of the
$\frak{Z}_n$ matrices and the properties of the preferred ordering on $P_n$.
\endproof

\bigskip
To evaluate the determinant $\det(\frak{Z}_n)$, we  follow the general advice of
\cite{K} and {\it introduce parameters} in the matrix entries. Specifically, we
consider the matrix:
$$\frak{Z}_n(x,y) = xZ_n + yZ_n^t = \pmatrix xZ_{n-1} + yZ_{n-1}^t &
xZ_{n-1}\cr
                                             yZ_{n-1}^t & xZ_{n-1} + yZ_{n-1}^t\cr\endpmatrix$$
Using the Maple computer algebra system, the determinants of the first few of these
matrices are found to be:
$$\matrix  n & \det(\frak{Z}_n(x,y)) \cr
& \cr 1 & x^2 + xy + y^2\cr 2 & (x+y)^2(x^2 - xy + y^2)\cr 3 &
(x-y)^2(x^2+xy+y^2)^3\cr 4 & (x+y)^6(x^2-xy+y^2)^5\cr 5 &
(x-y)^{10}(x^2+xy+y^2)^{11} \cr\endmatrix\leqno(4)$$ An interesting recurrence
explains the patterns evident in these examples.

\bigskip
\proclaim{Lemma 4}  The determinants of the $\frak{Z}_n(x,y)$ matrices are related
by the following recurrence:
$$\det(\frak{Z}_{n+2}(x,y)) = (\det(\frak{Z}_n(x,y))^2
\det(\frak{Z}_{n+1}(x,-y)).$$
(Note that the right side involves both
$\frak{Z}_{n+1}$ and $\frak{Z}_{n}$
and that $y$ is negated in the second
factor.)
\endproclaim

\bigskip
\noindent \proof  We use part $(2)$ of Lemma 3 twice to form the following block
decomposition of $\frak{Z}_{n+2}(x,y)$:
$$\frak{Z}_{n+2}(x,y) =
\pmatrix xZ_n + yZ_n^t & xZ_n & xZ_n & xZ_n\cr
                                  yZ_n^t & xZ_n +yZ_n^t & 0 & xZ_n\cr
                                  yZ_n^t & 0 & xZ_n +yZ_n^t & xZ_n\cr
                                  yZ_n^t & yZ_n^t & yZ_n^t & xZ_n +yZ_n^t\cr \endpmatrix$$
(where each entry is a block of size $2^n \times 2^n$).  To evaluate the
determinant,
we perform block-wise row and column operations.  To simplify
the notation,
we write $Z = Z_n$ First subtract
row 4 from each of the first three rows
to obtain:
$$\pmatrix xZ & xZ-yZ^t & xZ-yZ^t & -yZ^t\cr
            0 & xZ & -yZ^t & -yZ^t\cr
            0 & -yZ^t & xZ & -yZ^t\cr
         yZ^t & yZ^t & yZ^t & xZ +yZ^t\cr\endpmatrix.$$
Subtract column 1 from each of the columns 2,3,4 to obtain:
$$\pmatrix xZ & -yZ^t & -yZ^t & -xZ-yZ^t\cr
            0 & xZ & -yZ^t & -yZ^t\cr
            0 & -yZ^t & xZ & -yZ^t\cr
            yZ^t & 0 & 0 & xZ \cr\endpmatrix.$$
Next, subtract column 2 from column 4:
$$\pmatrix xZ & -yZ^t & -yZ^t & -xZ\cr
            0 & xZ & -yZ^t & -xZ-yZ^t\cr
            0 & -yZ^t & xZ & 0\cr
            yZ^t & 0 & 0 & xZ \cr\endpmatrix.$$
Add column 1 to column 4, then add row 4 in the resulting matrix to row 2:
$$\frak{Z}' = \pmatrix xZ & -yZ^t & -yZ^t & 0\cr
                        yZ^t & xZ & -yZ^t & 0\cr
                       0 & -yZ^t & xZ & 0\cr
                       yZ^t & 0 & 0 & xZ+yZ^t \cr\endpmatrix.$$
Expanding along the last column we obtain
$$\det(\frak{Z}_{n+2}(x,y)) = \det(\frak{Z}')
= \det(xZ+yZ^t)\det\pmatrix xZ &
 -yZ^t & -yZ^t\cr
                             yZ^t & xZ & -yZ^t\cr
                            0 & -yZ^t &
xZ\cr\endpmatrix.\leqno(5)$$ The first factor on the right of (5) is
$\det(\frak{Z}_n(x,y))$.
 Continuing with the $3\times 3$ matrix,
subtract column 3 from
column 2:
$$\pmatrix xZ & 0 & -yZ^t\cr
         yZ^t & xZ+yZ^t & -yZ^t\cr
          0 & -xZ-yZ^t & xZ\cr\endpmatrix,$$
then add row three to row 2:
$$\pmatrix xZ & 0 & -yZ^t\cr
          yZ^t & 0 & xZ-yZ^t\cr
          0 & -xZ-yZ^t & xZ\cr\endpmatrix.$$
Expanding along column 2, we have
$$\eqalign{\det(\frak{Z}_{n+2}(x,y)) &=
\det(\frak{Z}_n(x,y))^2
  \det\pmatrix xZ & -yZ^t\cr
                        yZ^t & xZ-yZ^t\cr\endpmatrix\cr
  &=\det(\frak{Z}_n(x,y))^2\det\pmatrix xZ-yZ^t & -yZ^t\cr
                             xZ & xZ-yZ^t\cr\endpmatrix \cr
  &=\det(\frak{Z}_n(x,y))^2\det(\frak{Z}_{n+1}(x,-y)),\cr}$$
as claimed.  For the last equality, we perform
row and column interchanges
to put the final matrix shown into the form:
$$\pmatrix xZ-yZ^t & xZ\cr
         -yZ^t & xZ-yZ^t\cr\endpmatrix$$
required for $\frak{Z}_{n+1}(x,-y)$.
\endproof

\bigskip
From the initial cases computed with Maple in (4) and the recurrence from Lemma 4,
we see that there are nonnegative integers $\alpha_n,\beta_n$ such that
$$\det(\frak{Z}_n(x,y)) = (x+(-1)^n y)^{\alpha_n} (x^2-(-1)^n
xy+y^2)^{\beta_n}.\leqno(6)$$ Moreover, the recurrence from Lemma 4 implies that
$$\cases \alpha_{n+2} = 2\alpha_n + \alpha_{n+1},&\cr
          \beta_{n+2} = 2\beta_n + \beta_{n+1}.&\cr\endcases\leqno(7)$$

We also have the following fact that is evident from (4):

\proclaim{Lemma 5}  For all $n\ge 1$, $\alpha_n = \beta_n + (-1)^{n+1}$.\endproclaim

\bigskip
\noindent \proof  This follows by induction.  The base cases come from the Maple
computations in (4) above: $\alpha_1 = 0, \beta_1 = 1$, and $\alpha_2 = 2, \beta_2 =
1$. For the induction step, assume that the claim of the lemma has been proved for
all $\ell \le k + 1$.  Then subtracting the two recurrences from (7) shows that
$$\alpha_{k+2} - \beta_{k+2} =
2(\alpha_k - \beta_k) + \alpha_{k+1} - \beta_{k+1} = 2(-1)^{k+1} + (-1)^{k+2} =
(-1)^{k+1}=(-1)^{k+3}.\ \ \square$$

\bigskip
We are now ready to conclude the proof of our main theorem. To determine the
determinant of the original $\frak{Z}_n = \frak{Z}_n(1,1)$, we simply substitute $x
= y = 1$ in  (6).  The factors of $x - y$ show immediately that $\det(\frak{Z}_n) =
0$ if $n$ is odd.  Moreover, when $n$ is even we have $\det(\frak{Z}_n) =
2^{\alpha_n}$.  We  solve the first recurrence in (7) for $\alpha_n$ by the standard
method for second order linear recurrences with constant coefficients.  The
characteristic equation is $r^2 - r - 2 = 0$, whose roots are $r = 2, -1$.  Hence
$\alpha_n = c_1(2)^n + c_2 (-1)^n$ for some constants $c_1,c_2$. The initial
conditions $\alpha_1=0, \alpha_2 = 2$ show that $c_1 = 1/3, c_2 = 2/3$.  Hence:
$$\alpha_n = {1\over 3}(2)^n + {2\over 3} (-1)^n.$$
Hence if $n$, and therefore also $n + 2$, are even, we have
$$\alpha_{n+2} = {1\over 3}(2)^{n+2} + {2\over 3} =
4\left({1\over 3}(2)^n + {2\over 3}\right) - 2
= 4\alpha_n - 2.\quad \square $$
\bigskip

\proclaim{Corollary 2} If $n$ is even, then
$$ \det({\frak{Z}_n})=2^{\textstyle{\frac{2^n+2}{3}}}.$$
\endproclaim
\bigskip

\refstyle{A}

\enddocument